\documentclass{amsart}
\usepackage{palatino}
\usepackage{amsthm, amsmath}
\usepackage{amssymb} 
\usepackage[margin=3cm]{geometry}
\usepackage[breaklinks,bookmarksopen,bookmarksnumbered]{hyperref}
\usepackage[all]{xy}

\theoremstyle{plain}
\newtheorem*{thm*}{Theorem}
\newtheorem{lem}{Lemma}
\newtheorem*{prop}{Proposition}
\theoremstyle{remark}
\newtheorem*{rem}{Remark}

\newcommand\pr{\noindent\textit{Proof} : }

\newcommand\Aut{\operatorname{Aut}}

\renewcommand\ss{\mathfrak{S}_6}

\newcommand\Z{\mathbb{Z}}

\newcommand\C{\mathbb{C}}
\renewcommand\P{\mathbb{P}}

\renewcommand\O{\mathcal{O}}

\renewcommand\ss{\mathfrak{S}_6}

\newcommand\iso{\vbox{\hbox to .8cm{\hfill{$\scriptstyle\sim$}\hfill}
\nointerlineskip\hbox to .8cm{{\hfill$\longrightarrow $\hfill}} }}

\begin{document}
\title[Non-rationality of the $\mathfrak{S}_6$-symmetric quartic threefolds]{Non-rationality of the $\mathfrak{S}_6$-symmetric quartic threefolds}
\author[Arnaud Beauville]{Arnaud Beauville}
\address{Laboratoire J.-A. Dieudonn\'e\\
UMR 7351 du CNRS\\
Universit\'e de Nice\\
Parc Valrose\\
F-06108 Nice cedex 2, France}
\email{arnaud.beauville@unice.fr}

\begin{abstract}
We prove that the quartic hypersurfaces defined by $\  \sum x_i=t\sum x_i^4-(\sum x_i^2)^2=0 $ in $\P^5 $ are not rational 
for $t\neq 0,2,4,6 ,\frac{10}{7}\,\cdot $
\end{abstract}
\maketitle 
\section{Introduction}
Let  $V$ be the standard representation of $\ss$ (that is, $V$ is the hyperplane $\sum x_i=0$ in $\C^6$, with $\ss$ acting by permutation of the basis vectors).
The quartic hypersurfaces in $\P(V)\ (\cong \P^4)$ invariant under $\ss$ form the pencil
\[X_t:\    t\sum x_i^4-(\sum x_i^2)^2=0 \ ,\quad t\in \P^1\ .\]
This pencil contains the Burkhardt quartic (for $t=2 $) and  the Igusa quartic ($t={4} $), which are both rational.

For  $t\neq 0,2,4,6 $ and $\frac{10}{17}$, the quartic $X_t$ has exactly 30 nodes; the set of nodes $\mathcal{N}$ is the orbit under $\ss$ of $(1,1,\rho ,\rho ,\rho ^2,\rho ^2)$, with $\rho =e^{\frac{2\pi i}{3} }$ (\cite{vdG}, \S 4). 
We will prove:
\begin{thm*}
For $t\neq 0,2,4,6 ,\frac{10}{7}$, $X_t$ is not rational.
\end{thm*}
The method is that of \cite{B} : we show that the intermediate Jacobian of a desingularization of $X_t$ is 5-dimensional and that the action of $\ss$ on its tangent space at $0$ is irreducible. From this one sees easily that this intermediate Jacobian cannot be a Jacobian or a product of Jacobians, hence  $X_t$ is not rational by the Clemens-Griffiths criterion. We do not know whether $X_t$ is unirational.
 
 \medskip	
{\small I am indebted to A. Bondal and Y. Prokhorov for suggesting the problem, and to A. Dimca for explaining to me how to compute explicitly the defect of a nodal hypersurface.}

\bigskip	
\section{The action of $\ss$ on $T_0(JX)$}

We fix $t\neq 0,2,4,6 ,\frac{10}{7}$, and denote by $X$ the desingularization of $X_t$ obtained by blowing up the nodes. The main ingredient of the proof is the fact that the action of $\ss$ on $JX$ is non-trivial. To prove this we consider the action of $\ss$ on the tangent space $T_0(JX)$, which is by definition $H^2(X,\Omega ^1_X)$.

\begin{lem}
Let $\mathcal{C}$ be the space of cubic forms on $\P(V)$ vanishing along $\mathcal{N}$. We have an isomorphism of $\ss$-modules $\mathcal{C}\cong V\oplus H^2(X,\Omega ^1_X)$.
\end{lem}

\pr The proof is essentially contained in \cite{C}; we explain  how to adapt the arguments there to our situation. Let $b: P\rightarrow \P(V)$ be the blowing-up of $\P(V)$ along $\mathcal{N}$. The threefold $X$ is the strict transform of $X_t$ in $P$.
The exact sequence 
\[0\rightarrow N^*_{X/P}\longrightarrow \Omega ^1_{P\, | X}\longrightarrow \Omega ^1_X\rightarrow 0 \]
gives rise to an exact sequence
\[0\rightarrow H^2(X,\Omega ^1_X)\longrightarrow H^3(X,N^*_{X/P}) \longrightarrow H^3(X,\Omega ^1_{P\, | X})\rightarrow 0\ . \] (\cite{C}, proof of theorem 1), which is $\ss$-equivariant. We will compute the two last terms. 

The exact sequence
\[0\rightarrow \Omega ^1_P(-X)\longrightarrow \Omega ^1_P\longrightarrow \Omega ^1_{P\, | X}\rightarrow 0\]
provides an isomorphism $H^3(X,\Omega ^1_{P\, | X})\iso H^4(P, \Omega ^1_P(-X))$, and the latter space is isomorphic to
$H^4(\P(V), \Omega ^1_{\P(V)}(-4))$ (\cite{C}, proof of Lemma 3). By Serre duality $H^4(\P(V), \Omega ^1_{\P(V)}(-4))$ is dual to $H^0(\P(V), T_{\P(V)}(-1))\cong V$. Thus the  $\ss$-module $H^3(X,\Omega ^1_{P\, | X})$ is isomorphic to $V^*$,  hence also to $V$.

Similarly the exact sequence $\ 0\rightarrow \O_P(-2X)\longrightarrow \O_P(-X)\longrightarrow N^*_{X/P}\rightarrow 0\ $ and the vanishing of  $H^i(P, \O_P(-X))$ (\cite{C}, Corollary 2) provide an isomorphism of $H^3(X,N^*_{X/P})$ onto $H^4(P, \O_P(-2X))$, which is  naturally isomorphic to the dual of $\mathcal{C}$ (\cite{C}, proof of Proposition 2). The lemma follows.\qed

\begin{lem}
The dimension of $\mathcal{C}$ is $10$. 
\end{lem}
\pr Recall that the \emph{defect} of $X_t$ is the difference between the dimension of $\mathcal{C}$ and its expected dimension, namely : 
\[ \mathrm{def}(X_t):= \dim \mathcal{C} -( \dim H^0(\P(V), \O_{\P(V)}(3) )-\#\,\mathcal{N})\ .\]
Thus our assertion is equivalent to $\mathrm{def}(X_t)=5$. 

To compute this defect we use the formula of \cite{D-S}, Theorem 1.5. Let $F=0$ be an equation of $X_t$ in $\P^4$; let $R:= \C[X_0,\ldots ,X_4]/(F'_{X_0},\ldots ,F'_{X_4})$ be the Jacobian ring of $F$, and let $R^{sm}$ be the Jacobian ring of a \emph{smooth} quartic hypersurface  in $\P^4$. The formula is
\[ \mathrm{def}(X_t)= \dim R_7-\dim R^{sm}_7\ .\]
In our case we have $\dim R^{sm}_7=\dim R^{sm}_3=35-5=30$; a simple computation with Singular (for instance) gives $\dim R_7=35$. This implies the lemma.\qed
\bigskip	
\begin{prop}
The $\ss$-module $H^2(X,\Omega ^1_X)$ is isomorphic to $V$.
\end{prop}

\pr 
Consider the homomorphisms $a$ and $b$  of $\C^6$ into $H^0(\P(V), \O_{\P(V)}(3))$ given by $a(e_i)=x_i^3$, $b(e_i)=x_i\sum x_j^2$. They are both $\ss$-equivariant and map $V$ into $\mathcal{C}$;  the subspaces $a(V)$ and $b(V)$ of $\mathcal{C}$ do not coincide, so we have $a(V)\cap b(V)=0$. By Lemma 2 this implies $\mathcal{C}=a(V)\oplus b(V)$, so 
$H^2(X,\Omega ^1_X)$ is isomorphic to $V$ by Lemma 1.\qed
\medskip	
\begin{rem}
Suppose $t=2,6$ or $\frac{10}{7} $. Then the singular locus of $X_t$ is $\mathcal{N}\cup\mathcal{N}'$, where $\mathcal{N}'$ is the $\ss$-orbit of the point $(1,-1,0,0,0,0)$ for $t=2$, $(-1,-1,-1,1,1,1)$ for $t=6$,  $(-5,1,1,1,1,1)$ for $t=\frac{10}{7} $ \cite{vdG}. Since $x_1^3-x_0^3$ does not vanish on $\mathcal{N}'$, the space  of cubics vanishing along $\mathcal{N}\cup\mathcal{N}'$ is strictly contained in $\mathcal{C}$. 
By Lemma 1 it contains a copy of $V$, hence it is isomorphic to $V$; therefore  $H^2(X,\Omega ^1_X)$ and $JX$ are zero in these cases. We have already mentioned that $X_2$ and $X_4$ are rational;  we do not know whether this is the case for $X_6$ and $X_{\frac{10}{7} }$.
 
\end{rem}

\bigskip	
\section{Proof of the theorem}
To prove that $X$ is not rational, we apply  the Clemens-Griffiths criterion (\cite{CG}, Cor. 3.26): it suffices to prove that $JX$ is not a Jacobian or a product of Jacobians.

Suppose $JX\cong JC$ for some curve $C$ of genus $5$. By the Proposition $\ss$ embeds into the group  of automorphisms of $JC$ preserving the principal polarization; by the Torelli theorem this group is isomorphic to $\Aut(C)$ if $C$ is hyperelliptic and $\Aut(C)\,\times \,\Z/2\ $ otherwise. Thus we find $\#\Aut(C)\geq \frac{1}{2}6! =360$. But this contradicts the Hurwitz  bound $\ \#\Aut(C)\leq 84(5-1)=336$.

Now suppose that $JX$ is isomorphic to a product of Jacobians $J_1\times \ldots \times J_p$, with $p\geq 2$.  Recall that such a decomposition is \emph{unique} up to the order of the factors: it
corresponds to the decomposition of the Theta divisor into irreducible components (\cite{CG}, Cor. 3.23). Thus the group $\ss$ permutes the factors $J_i$, and therefore acts on $[1,p]$; by the Proposition this action must be transitive. But we have $p\leq \dim JX=5$, so this is impossible.\qed

\end{document}